
\input amstex
\documentstyle{amsppt}
\def\IQ{\bold Q}
\def\IZ{\bold Z}
\def\IR{\bold R}
\def\IN{\bold N}
\def\a{\alpha}
\def\e{\varepsilon}
\def\ti{\times}
\def\id{\operatorname{id}}
\def\Iso{\operatorname{Iso}}
\def\Auth{\operatorname{Auth}}
\def\bs{\setminus}
\NoBlackBoxes

\topmatter
\title
Symmetry and colorings:\\
Some results and open problems
\endtitle
\author
Taras BANAKH and Igor PROTASOV
\endauthor
\subjclass 05D10
\endsubjclass
\keywords
coloring, symmetry, measure, cardinality, monochromatic subset,
symmetric subset, group, geometric figure
\endkeywords
\address Department of Mathematics, Lviv University of I.Franko, Universytetska 1, Lviv,
Ukraine\endaddress
\email tbanakh\@franko.lviv.ua
\endemail
\address Faculty of Cybernetics, Kiev University of T.Shevchenko, Volodymyrska 64, Kiev, 
Ukraine 
\endaddress
\email igor\@protasov.kiev.ua\endemail
\abstract We survey some principal results and open problems related to
colorings of algebraic and geometric objects endowed with symmetries.
\endabstract
\endtopmatter

\document

In this paper we give a survey of principal results and open problems related
to colorings and symmetry. 
Starting in 1995 from seminal questions of I.V.Protasov
\cite{13},\cite{15}, the subject
quickly tranformed into a field which can be considered as an independent 
branch of Ramsey theory with its specific methods, non-trivial problems and 
deep ties with many other mathematical disciplines such as
combinatorial set theory, theory of (compact) groups, general and algebraic 
topology, elementary and Riemannian geometry, probability theory, theory
of measure, functional and abstract harmonic analysis.

Roughly speaking, a typical problem in the subject sounds as follows: given a
space $X$ endowed with symmetries, find the maximal size of a monochromatic
symmetric subset $A$ which can be found in each ``good'' coloring of $X$. The
size of $A$ can be understood twofold: in the sense of cardinality as well as
of measure. Respectively, we have two groups of results coinciding on finite
objects. It turns out that non-trivial questions appear even for
simplest objects such as the unit interval $[0,1]$.

In the paper we survey results obtained in the last 5 years,
see \cite{1}--\cite{18}. 
We would like to express our gratitude to all mathematicians
working in the field. Especially we appreciate the contribution of
Ya.Vorobets and O.Verbitsky who have enriched the theory by many interesting 
and unexpected results.
In spite of some achievements, the field still contains many interesting 
unsolved problems and we hope that this survey will accelerate 
their solution.

\heading
1. Big monochromatic symmetric subsets in colorings of groups.
\endheading

Working with resolvable groups (see the survey \cite{14} for more
information on this subject), I.V.Protasov came in 1995 to the 
following question (see \cite{13} and \cite{15}): 
{\it is it true that for every $n$-coloring
of the group $\IZ^n$ there exists an infinite monochromatic subset of $\IZ^n$, 
symmetric with respect to some point of $\IZ^n$?} 
Elementary arguments show that the answer
is positive for $n=2$. For higher $n$ the problem turned to be not elementary
and for its solution required involving the machinery of algebraic topology
(Borsuk-Ulam Theorem and the theory of degree of maps between spheres).

We reformulate the Protasov question in a more general setting. Note that a
subset $A\subset \IR^n$ is symmetric with respect to a point $x\in\IR^n$ iff
$A=2x-A$. Keeping this in mind we define a subset $A$ of a group $G$ to be
{\it symmetric} if $A=gA^{-1}g$ for some $g\in G$. For a group $G$ let
$\nu(G)$ be {\it the minimal number of colors necessary to color the group
$G$ so that $G$ contains no infinite monochromatic symmetric subset}.
Clearly, the Protasov question is equivalent to asking if $\nu(\IZ^n)>n$.

Abelian groups $G$ with $\nu(G)=2$ were characterized in \cite{13}. In
\cite{10} the invariant $\nu(G)$ was calculated for any Abelian group $G$. It
can be expressed via known group invariants such as the free rank $r_0(G)$,
the cardinality of the group $G$, and the cardinality of its subgroup
$B(G)=\{g\in G: 2g=0\}$ consisting of elements of order 2.

\proclaim{1.1. Theorem (\cite{10})} For any Abelian group $G$
$$
\nu(G)=\cases
r_0(G)+1,&\text{if $G$ is finitely generated};\\
r_0(G)+2,&\text{if $G$ is countable, infinitely generated, 
and $|B(G)|<\aleph_0$};\\
\max\{|B(G)|,\log|G|\},&\text{if $G$ is uncountable or $|B(G)|\ge\aleph_0$};
\endcases
$$
where $\log|G|=\min\{\alpha:2^\alpha\ge|G|\}$.
\endproclaim

Thus Theorem 1.1 completely resolves the problem of calculating the
invariant $\nu(G)$ for Abelian groups. In particular, we get
$\nu(\IZ^n)=n+1$, $\nu(\IQ^n)=n+2$ and $\nu(\IR^n)=\aleph_0$ for every
$n\in\IN$. 

\proclaim{1.2. Problem ({\rm R.I.Grigorchuk})} Investigate the invariant
$\nu(G)$ for non-commutative groups. In particular, is $\nu(F_2)$ finite,
where $F_2$ is the free group with two generators?
\endproclaim

The only information on $\nu(F_2)$ we know is that $\nu(F_2)>2$, see \cite{GK}. 
What concerns big {\it finite} symmetric monochromatic subsets, we have the
following result due to I.V.Protasov.

\proclaim{1.3. Theorem (\cite{13})} For every finite coloring of an
infinite Abelian group $G$ there exists a symmetric monochromatic subset
$A\subset G$ of arbitrary big finite size.
\endproclaim

In fact, a more precise statement is true, see Theorem 5.1.

\proclaim{1.4. Theorem} If $r,n$ are natural numbers, then for every
$r$-coloring of an Abelian group $G$ with $|G|\ge r^2n$ there exists a
symmetric monochromatic subset $A\subset G$ of size $|A|\ge n$.
\endproclaim

A similar result holds also for finite non-commutative groups. 
Below $\pi(x)$ is the function assigning to each
$n$ the number of prime numbers not exceeding $n$.

\proclaim{1.5. Theorem} If $r,n$ are natural numbers, then for every
$r$-coloring of a finite group $G$ with $|G|>r^2(n!)^{\pi(n)}$
there exists a symmetric monochromatic subset $A\subset G$ of size
$|A|>n$. 
\endproclaim

This theorem follows from Theorem 1.4 and the known fact stating that each
finite group $G$ with $|G|>(n!)^{\pi(n)}$ contains an Abelian
subgroup $H$ of size $|H|>n$. In its turn, this fact can be easily
derived from the Silov Theorem 11.1 and Theorem 16.2.6 of \cite{23}. Note
that in contrast to Theorem 1.4, Theorem 1.5 says nothing about infinite
groups. 

\proclaim{1.6. Problem ({\rm I.Protasov})} Does for every finite coloring
of an infinite group $G$ there exist a monochromatic symmetric subset
$A\subset G$ of arbitrary big finite size?
\endproclaim

Note that this problem does not reduce to the Abelian case since there
exist infinite periodic groups containing no big Abelian subgroups (for
example, for every odd $d\ge 665$ and every $m\ge 2$ the free Bernside
group $B(m,d)$ is infinite and contains no Abelian subgroup of size $>d$,
see \cite{19}). In the meantime, by the theorem of Kargapolov-Hall-Kulatilani
(see \cite{24}), each locally finite infinite group does contain an
infinite Abelian subgroup and thus for locally finite groups Problem 1.6
has positive answer.

\heading
2. Central sets in Abelian groups
\endheading

Answering the Protasov question \cite{15}, T.Banakh noticed that for every
$n$-coloring of the group $\IZ^n$ there exists an infinite monochromatic
subset $A\subset \IZ^n$ symmetric with respect to some point $x\in\IZ^n$
of the cube $[-1,1]^n\subset\IR^n$. Having this in mind, we define a
subset $E$ of a group $G$ to be {\it central}\footnote{In [2] central sets are referred to as essential sets} if for every $r$-coloring of
$G$, where $r<\nu(G)$, there exists an infinite monochromatic subset
$A\subset G$, symmetric with respect to some point $x\in E$, i.e., 
$A=xA^{-1}x$. A central set
$E$ of $G$ is called a {\it minimal central set} in $G$ if no proper
subset of $E$ is a central set in $G$.

In this section given a group $G$ we consider two questions:
(1) what is the structure of central subsets of $G$ and (2) what is the
minimal size $c(G)$ of a central set in $G$?
While the first question has
appeared to be rather difficult (a satisfactory answer is known only for the
groups $\Bbb Z$, $\Bbb Z^2$, and $\Bbb Z^3$), the second question admits the
following answer (below $\frak c$ denotes the cardinality of continuum).

\proclaim{2.1. Theorem (\cite{2})}
For an infinite Abelian group $G$ of cardinality $\le\frak c$
$$
c(G)=\cases
1, &\text{if $|B(G)|\ge\aleph_0$ or $\nu(G)=2$};\\
\chi(\Bbb Z^{r_0(G)}),&\text{if $G$ is finitely generated};\\
\aleph_0,& \text{otherwise}.
\endcases
$$
\endproclaim

Thus, the problem of calculating the cardinal $c(G)$ reduces to
calculating the numbers $c(\Bbb Z^n)$ for the groups $\Bbb Z^n$.
In this respect, we
possess the following (unfortunately incomplete) information.

\proclaim{2.2. Theorem (\cite{2})} For every $n\in\Bbb N$
\roster
\item $c(\Bbb Z^n)\le c(\Bbb Z^{n+1})$;
\item $\frac{n(n+1)}2\le c(\Bbb Z^n)<2^n$;
\item $c(\Bbb Z^n)=\frac{n(n+1)}2$ if $n\le 3$.
\endroster
\endproclaim

Thus, only for $n\le 3$ the exact values of the cardinals $c(\Bbb Z^n)$
are known. These values are calculated by studying the structure of
minimal central subsets of the groups $\Bbb Z^n$ which are considered as
subgroups of linear spaces $\Bbb R^n$. Below under a {\it triangle} in $\Bbb
R^2$ we understand any 3-element affinely independent subset of $\Bbb
R^2$; an {\it octahedron} in $\Bbb R^3$ is any 6-element subset of the
form $\{a\pm e_i:i=1,2,3\}$ for some point $a\in\Bbb R^3$ and some
linearly independent vectors $e_1, e_2, e_3$ in $\Bbb R^3$. The following
statement characterizes triangles and octahedra as minimal central
subsets in $\Bbb Z^2$ and $\Bbb Z^3$.

\proclaim{2.3. Theorem (\cite{2})}
\roster
\item Triangles and only triangles are minimal central subsets of the
group $\Bbb Z^2$;
\item Every central subset of $\Bbb Z^2$ contains a triangle;
\item Octahedra and only octahedra are central subsets of
cardinality $\le 6=\chi(\Bbb Z^3)$ in the group $\Bbb Z^3$;
\item There is a finite central subset in $\Bbb Z^3$ containing no
octahedron. 
\endroster
\endproclaim

This theorem implies that for $n=2,3$ central subsets of cardinality
$\chi(\Bbb Z^n)$ in $\Bbb Z^n$ can be realized as the sets of vertices of
respective regular polyhedra in $\Bbb R^n$. This already is not true for
$n=4$: no regular polyhedron in $\Bbb R^4$ exists such that its set of
vertices is an central subset in $\Bbb Z^4$ of cardinality $\chi(\Bbb
Z^4)$. This is so because $11\le\chi(\Bbb Z^4)\le 14$ (see \cite{2}),
while regular polyhedra in $\Bbb R^4$ contain 5, 8, 16, 24, 120, or 600
vertices, see \cite{20, 12.6}.

\proclaim{2.4. Problem (\cite{2})} Calculate the numbers $\chi(\Bbb Z^n)$ and
investigate  the structure of central sets in $\Bbb Z^n$ for $n>3$.
\endproclaim

\heading
3. Uncountable monochromatic symmetric subsets in colorings of Abelian
groups 
\endheading

In this section we deal with the following modification of the question
considered in the first section. {\it Let $G$ be a group and $\kappa,r$ be
cardinal numbers. Does for every $r$-coloring of $G$ there exist a
monochromatic symmetric subset of cardinality $\kappa$?}

First we state the following result of I.V.Protasov.

\proclaim{3.1. Theorem (\cite{16})} For every uncountable Abelian group
$G$ with $|B(G)|<|G|$ there exists a 2-coloring of $G$ having no
monochromatic symmetric subset of cardinality $|G|$.
\endproclaim

It will be convenient for every group $G$ and a cardinal $\kappa$ to introduce
the following cardinal invariant generalizing $\nu(G)$: let $\nu(G,\kappa)$ be
the minimal cardinal number of colors necessary to colorate the group $G$ so
that $G$ would contain no monochromatic symmetric subset of cardinality
$\ge \kappa$. Clearly, $\nu(G)=\nu(G,\aleph_0)$ and Theorem 3.1 states that
$\nu(G,|G|)=2$ for every uncountable Abelian group $G$ with $|B(G)|<|G|$.

The problem of calculating the cardinals $\nu(G,\kappa)$ turned out to be 
tightly connected with some well-known partition relations for cardinals,
see \cite{21}. For a set $S$ by $[S]^2$ the collection of all 2-element subsets
of $S$ is denoted. For ordinals $\lambda,\a$ and a cardinal $r$ the symbol
$\lambda\to(\a)^2_r$ means that for any $r$-coloring the set $[\lambda]^2$
there exists a subset $A\subset\lambda$ of order type $\a$ such that the
set $[A]^2\subset[\lambda]^2$ is monochromatic. We identify the cardinals with
the smallest ordinals of the corresponding size.

The following theorem generalizes a result of \cite{16} and 
can be proved by analogy
with the proof of Lemmas 13--15 of \cite{10}.

\proclaim{3.2. Theorem} If $\lambda,\kappa,r$ are cardinals such that
$\lambda\to(\kappa+1)^2_r$, then $\nu(G,\kappa)>r$ for every infinite
Abelian group $G$ of size $|G|\ge\lambda$, that is for every $r$-coloring
of the group $G$ there exists a monochromatic
symmetric subset $A\subset G$ of size $|A|\ge k$.
\endproclaim

It follows from the famous Erd\"os-Rado stepping up lemma \cite{21}
that $(2^{<\kappa})^+\to(k+1)^2_r$ for every cardinals
$r<\kappa\ge\aleph_0$, where $2^{<\kappa}=\sup\{2^\tau:\tau<\kappa\}$.
This observation and Theorem 3.2 imply

\proclaim{3.3. Corollary} $\nu(G,\kappa)\ge k$ for every cardinal $\kappa$
and every infinite Abelian group $G$ of size $|G|>2^{<\kappa}$.
\endproclaim

Observe that GCH, the Generalized Continuum Hypothesis, is equivalent to the
assumption $\kappa=2^{<\kappa}$ for every infinite cardinal $\kappa$.
Thus under GCH, $\nu(G,\kappa)\ge\aleph_0$ for every cardinal $\kappa$ and
every infinite Abelian group $G$ of size $|G|>\kappa$. We do not know if
this result is true in ZFC. The following result suggests a positive
answer to this question.

\proclaim{3.4. Proposition} For every infinite group $G$ and every
cardinal $\kappa<|G|$ we have $\nu(G,\kappa)>3$, that is for every
3-coloring of $G$ there exists a symmetric monochromatic subset $A\subset
G$ of size $|A|\ge \kappa$.
\endproclaim

This proposition  can be easily derived from the next modification of
Lemma 1 of \cite{16}.

\proclaim{3.5. Lemma} Let $G$ be an infinite Abelian group, $k<|G|$, and
$H\subset G$ be a subgroup with $|H|\le k$. For every $k$-coloring
$\chi:G\to k$ of the group $G$ one of the following statements holds:
\roster
\item $G$ contains a monochromatic symmetric subset of size $\ge k$;
\item there exists $g\in G$ such that
$\chi(2H+g)\cap\chi(2H-g)=\emptyset$. 
\endroster
\endproclaim

\heading
4. Invariant $ms(G,r)$ for finite groups.
\endheading

Observe that the invariant $\nu(G,k)$ is defined also for finite groups!
In this case, however, it is more convenient to consider the invariant
$MS(G,r)$ which is in some sense inverse to the invariant $\nu(G,r)$. For a group
$G$ and a number $r$ define $MS(G,r)$ to be the least upper bound of
cardinal numbers $\tau$ such that for every $r$-coloring of the group $G$
there exists a monochromatic symmetric subset $A\subset G$ of size
$|A|\ge\tau$. For finite groups $G$ it is more preferable to work with the
normed value $ms(G,r)=\frac{MS(G,r)}{|G|}$.

To estimate the invariant $ms(G,r)$, for every group $G$ we
introduce the following characteristics:
$$
k(G)=\min_{g\in G}|\sqrt{g^2}|,\quad k_0(G)=|\sqrt 1|,\quad\text{and}\quad 
k_1(G)=\min_{g\in 2G\bs\{1\}}|\sqrt g|,
$$
where 1 is the neutral element of the group $G$, $2G=\{g^2 :g\in
G\}\subset G$ and $\sqrt g=\{x\in G:x^2=g\}$ for $g\in G$. Let also $m(G)$
be the maximal size of a subgroup of $G$ contained in the set $2G$.

\proclaim{4.1. Theorem (\cite{4})} For every finite group $G$ and every
$r\in\IN$ we have the estimates:
\roster
\item
$ms(G,r)\le\frac1{r^2}+(\frac1r-\frac1{r^2})\frac{k_0(G)}{|G|}+3\sqrt{2\ln(2r|G|)/|G|}$;
\item $ms(G,r)\ge \frac{k(G)m(G)}{|G|}\cdot\frac1{r^2}$;
\item $ms(G,r)\ge
\frac{k_1(G)m(G)}{|G|}\frac1{r^2}+\frac{k_0(G)-k_1(G)}{|G|}\frac1{r}$.
\endroster
\endproclaim

For some groups (for example, dihedral groups) the third item
gives a better lower bound than the second one.  
For any Abelian group $G$ we have $m(G)=|2G|$ and
$k(G)=|B(G)|=\frac{|G|}{m(G)}$ and thus
$ms(G,r)\ge\frac1{r^2}$ (compare to Theorem 5.1). I.V.Protasov has
noticed that this inequality is true also for every group $G$ of odd order
since for such a group $G$ we have $m(G)=|G|$ and $k(G)=1$. Yet, for
Abelian groups of odd order we have more:

\proclaim{4.2. Theorem (\cite{7, \S2})} $ms(G,r)>\frac1{r^2}$ for every
Abelian group $G$ of odd order and every $r>1$ (that is for every
$r$-coloring of $G$ there exists a symmetric monochromatic subset of size
$>\frac1{r^2}|G|$).
\endproclaim

This theorem was proved by Ya.Vorobets by application of the Fourrier
transform. Since the Fourrier transform does not work properly in the
non-Abelian case, we come to

\proclaim{4.3. Question} Is $ms(G,r)>\frac1{r^2}$ for every group $G$ of
odd order and every $r>1$?
\endproclaim

Surprisingly, but for groups of even order the inequality
$ms(G,r)\ge\frac1{r^2}$ does not hold anymore: it was noticed in \cite{4}
that $ms(Q_8,2)=\frac18<\frac1{2^2}$, where $Q_8=\{\pm1,\pm i,\pm j,\pm
k\}$ is the group of quaternions. Nonetheless the following question
(related to Theorem 1.5) is still open.

\proclaim{4.4. Question (\cite{4})} Does for every $\e>0$ and $r>1$ there
exist a finite group $G$ with $ms(G,r)<\e$?
\endproclaim

\heading
5. Invariant $ms(G,r)$ for compact topological groups
\endheading

Theorem 4.1 implies that for every $r\in\IN$ the limit $\lim_{n\to\infty}
ms(\IZ_n,r)$, where $\IZ_n$ are cyclic groups of order $n$, exists and is
equal to $\frac1{r^2}$. The value $\frac1{r^2}$ has quite real sense. It
turns out that the invariant $ms(G,r)$ can be defined for every compact
topological group $G$ and $\frac1{r^2}$ is just equal to 
$ms(S^1,r)$, where $S^1=\IR/\IZ$ is the unit circle in the
complex plane. 

It is well known that each compact topological group $G$ admits a unique
Haar measure $\mu$ (that is a probability bi-invariant regular Borel
measure on $G$). For every such a group $G$ and every $r\in\IN$ let
$ms(G,r)$ be the least upper bound of $\e>0$ such that for every
measurable $r$-coloring of the group $G$ there exists a measurable
monochromatic symmetric subset $A\subset G$ of measure $\mu(A)\ge \e$. 

Obviously, this definition of $ms(G,r)$ agrees with that given in the
previous section for finite groups. Now it becomes clear what the symbol {\it
ms} means: it is the abbreviation for ``{\it m}aximal {\it m}easure of {\it
m}onochromatic {\it s}ymmetric {\it s}ub{\it s}et''.

Note that in the definition of $ms(G,r)$ we consider only {\it measurable}
coloring of $G$, that is
colorings whose monochromatic classes are measurable subsets of $G$ with
respect to the Haar measure $\mu$. (Admitting also non-measurable
colorings would yield $ms(S^1,r)=0$  according to Theorem 3.1). 

\proclaim{5.1. Theorem (\cite{7, 2.2})} For every compact Abelian group
$G$ and every $r\in\IN$ we have the following estimates:
$$
\frac1{r^2}\le ms(G,r)\le
\frac1{r^2}+\Big(\frac1r-\frac1{r^2}\Big)\mu(B(G)).
$$
\endproclaim

Note that $\mu(B(G))=0$ for every non-degenerate connected compact 
Abelian group $G$. 
Therefore, for such groups $G$ Theorem 5.1 gives the exact value 
$ms(G,r)=\frac1{r^2}$.

The problem of calculating the invariant $ms(G,r)$ for connected 
non-commutative groups is much more complex. We begin with connected Lie
groups. 

\proclaim{5.2. Theorem (\cite{5})} If $r\in\IN$ and $G$ is a compact
connected Lie group, then
$$
\frac1{r^2}\cdot \frac{k(G)}{2^{\dim(G)}}\le ms(G,r)\le\frac1{r^2},
\quad\text{where $k(G)=\min_{g\in G}|\sqrt{g^2}|$}.
$$
\endproclaim

It is interesting to note that $k(S_1^n)=2^n$ for the torus $S_1^n$
and thus Theorem 5.2 gives
the exact value for $ms(S_1^n)$! For arbitrary
connected compact group we have a weaker lower bound.

\proclaim{5.3. Theorem (\cite{5})} If $r\in\IN$ and $G$ is a
non-degenerate connected compact group, then
$$
\frac1{r^2}\cdot\frac1{2^{\dim G}}\le ms(G,r)\le\frac1{r^2}.
$$
\endproclaim

\proclaim{5.4. Question (\cite{5})} Is there a connected compact group
$G$ with $ms(G,r)<\frac1{r^2}$ for some $r\in\IN$? In particular, what
is $ms(SO(3),2)$ equal to?
\endproclaim

Note that Theorem 5.2 yields the estimate $\frac1{16}\le
ms(SO(3),2)\le\frac14$. 

\heading
6. Invariant $ms(X,S,r)$ for some geometric objects.
\endheading

Analyzing the definition of the invariant $ms(G,r)$ for compact groups, we
can see that it is defined in a much more general situation. Let $X$ be a
topological space equipped with a regular Borel measure $\mu$. By
$\Auth(X)$ we denote the group of all measure preserving homeomorphisms of
$X$. Subsets of $X$ measurable with respect to the measure $\mu$ will be
called {\it measurable subsets} of $X$. 

 Suppose we are given with a subset $\Cal S\subset\Auth (X)$ whose
elements are called {\it admissible symmetries} of $X$. We define a subset
$A\subset X$ to be {\it $\Cal S$-symmetric} if $A=s(A)$ for some $s\in
\Cal S$.

Given a measurable subset $A\subset X$ and a finite $r$ we define
$MS(A,\Cal S,r)$ be the least upper bound of $\e>0$ such that for every
measurable $r$-coloring of $A$ there exists a measurable monochromatic $\Cal
S$-symmetric subset $B\subset A$ of measure $\mu(B)\ge \e$. If
$0<\mu(A)<\infty$ we let $ms(A,\Cal S,r)=\frac{MS(A,\Cal S,r)}{\mu(A)}$.

Note that the invariant $ms$ is monotone in sense that $ms(A,\Cal S',r)\le
ms(A,\Cal S,r)$ for any $\Cal S'\subset\Cal S$.

Given a subset $X$ of a finite-dimensional Euclidean space let $E$ be 
the affine hull of $X$. There are many natural sets which can serve
as sets of admissible symmetries of $X$. The greatest is the group 
$\Iso(X)$ of all isometries $f$ of the space $E$ such that $f(X)=X$. 
All others are subsets in $\Iso(X)$: $\Iso^+(X)$ is the subgroup of $\Iso(X)$
consisting of all orientation-preserving isometries; $\Iso_2(X)$ is the set
of all involutive isometries of $X$, $\Cal I=\Iso(X)\bs \{\id\}$ is the set
of all nontrivial isometries of $X$, $\Cal I^+=\Cal I\cap \Iso^+(X)$, $\Cal
I^-=\Iso(X)\bs\Iso^+(X)$, $\Cal I_2=\Cal I\cap \Iso_2(X)$, 
$\Cal I_2^+=\Cal I_2\cap \Iso^+(X)$, and $\Cal
I_2^-=\Cal I_2\cap \Iso^-(X)$. 

It is known that each compact set $X$ admits a probability Borel measure
$\mu$ invariant with respect to all isometries, that is $\mu(f(A))=\mu(A)$
for any measurable subset $A\subset X$ and any isometry $f\in \Iso(X)$, see
\cite{22, 2.7.3}. We assume that any of considered below geometric figures
(including $n$-dimensional spheres $S^n$ and balls $B^n$) are
endowed with such an invariant measure. 

\proclaim{6.1. Theorem (\cite{7, \S3})} For every $n,r\in\IN$ we have
$$
\frac1{r^2}=ms(S^n,\Cal I,r)=ms(S^n,\Cal I_2^-,r)=ms(B^{n+1},\Cal
I,r)=ms(B^{n+1},\Cal I_2^-,r).
$$
\endproclaim

The proof of Theorem 6.1 involves the following general result related to
colorings of Riemannian manifolds. We note that each Riemannian manifold
$M$ possesses a canonical metric and measure (generated by the Riemannian
structure of $M$). 

\proclaim{6.2. Theorem (\cite{7, 4.1})} For every $r\in\IN$ and every
compact subset $X$ of positive measure in a connected Riemannian manifold
$M$ we have $ms(X,\Cal I,r)\le \frac1{r^2}$, where $\Cal I$ is the set of 
all non-identity isometries $f$ of $M$ such that $f(X)=X$.
\endproclaim

As noted in \cite{7, \S3} for every 
measurable $r$-coloring of the $n$-dimensional sphere $S^n$ there exists 
a measurable monochromatic subset $A\subset S^n$ of measure $>\frac
1{r^2}$, symmetric with respect to some hyperplane passing through the
center of the sphere.

\proclaim{6.3. Question} Is $ms(S^n,\Cal I_2^+,r)=\frac1{r^2}$ for 
$n,r>1$? More precisely, does for every measurable $r$-coloring of the
$n$-dimensional sphere $S^n$ there exist a monochromatic subset of measure
$\ge\frac1{r^2}$, symmetric with respect to some line passing through the
center of the sphere?
\endproclaim

It can be shown that $ms(S^2,\Cal I^+_2,r)\ge\frac1{2r^2}$ for every
$r\in\IN$ (use arguments analogous to those from \cite{5}). Moreover,
applying Theorems 7.1 and 7.2 from the next section 
we may conclude that $ms(S^2,\Cal I_2^+,2)\ge
ms(S^2,\Cal V,2)=\frac16$ and $ms(S^2,\Cal I_2^+,3)\ge ms(S^2,\Cal
V,3)\ge\frac1{18}$, where $\Cal V$ is a 3-element set consisting of
rotations on 180$^\circ$ around 3 pairwise orthogonal axes passing
through the center of the sphere $S^2$.

\heading
7. Invariant $ms(X,r)$ for geometric figures with finite isometry group.
\endheading

It turns out that the invariant $ms(X,\Cal I,r)$ is positive not only for
spheres and balls, but also for many geometric figures with finite
isometry group. We consider the problem from a bit more general point of view.

Let a finite group $G$ act on a compact space $X$ endowed with a
probability Borel measure $\mu$, invariant under this action.
For every finite $r$ we
define $ms(X,r)$ to be the least upper bound of $\e>0$ such that for every
measurable $r$-coloring of $X$ there exists a monochromatic subset
$A\subset X$ of measure $\mu(A)\ge\e$ such that $A=g\cdot A$ for some
non-identity element $g\in G$. 

Below $I$ denotes the interval $[0,1]$ endowed with the standard Lebesgue
measure $\lambda$. For every compact group $G$ we consider the cylinder
$G\ti I$ as a left $G$-space (with the action $g\cdot(h,t)=(gh,t)$ for $g\in
G$ and $(h,t)\in G\ti I$) endowed with the product measure $\mu\times
\lambda$, where $\mu$ is the Haar measure on $G$.

\proclaim{7.1. Theorem} Suppose $G$ is a finite group, $X$ is a 
metrizable compact $G$-space, and $\mu$ is a probability $G$-invariant
Borel measure on $X$. For every finite $r$ we have
\roster
\item $ms(X,r)\ge ms(G\ti I,r)$;
\item $ms(X,r)=ms(G\ti I,r)$ if the measure $\mu$ is continuous and
$|G\cdot x|=|G|$ for almost all $x\in X$.
\endroster
\endproclaim

We recall that a measure $\mu$ on $X$ is called {\it continuous} if
$\mu(\{x\})=0$ for every $x\in X$. Thus the problem reduces to calculating
the invariant $ms(G\ti I,r)$ for finite groups $G$. Below $\Cal D_{2n}$ is
a dihedral group, i.e., the group of all isometries of the regular
$n$-gon. For $n=2$ \ $\Cal D_{2n}$ coincides with the Kleinian group 
(isomorphic to $\IZ_2\ti\IZ_2$).

The following unexpected result belongs to Ya.Vorobets.

\proclaim{7.2. Theorem}
\roster
\item For any subgroup $H$ of a finite group $G$ we have $ms(H\ti I,r)\le
ms(G\ti I,r)$ for $r\in\IN$.
\item $ms(\IZ_n\ti I,r)=0$ for every $n,r>1$.
\item For every prime number $p$ \ $ms(\Cal D_{2p}\ti
I,2)=\frac{p-1}{p^2+2p-2}$, $ms(\Cal D_{2p}\ti I,3)=\frac1{3p^2+6}$ and
$ms(\Cal D_{2p},r)=0$ if $r\ge 4$.
\item For any $n,r>1$ \ $ms(\Cal D_{2n},r)=ms(\Cal D_{2p},r)$, where $p$
is the minimal prime number dividing $n$;
\item If $r=2^k$ for some $k\ge 0$, then
$ms(\IZ_2^n\ti I,r)=\frac1{r^2}\cdot\frac{2^n-r}{2^n-1}$ for every $n>k$.
\endroster
\endproclaim

Theorems 7.1 and 7.2 allow us to classify the invariant 
$ms(X,r)$ for all compact
subsets of positive Lebesgue measure in the plane, where $ms$ is related
to the set $\Cal I=\Iso(X)\bs\{\id\}$ of all nonidentity isometries of $X$.

\proclaim{7.3. Theorem} For every compact subset $X\subset\IR^2$ of positive 
Lebesgue measure  and  every $r\ge 2$ we have
$$
ms(X,r)=\cases
\frac1{r^2},&\text{if the group \ $\Iso(X)$ is infinite;}\\
ms(\Cal D_{2n}\ti I,r),&\text{if \ $\Iso(X)$ is isomorphic to}\\
&\text{the dihedral group $\Cal D_{2n}$ for some $n\ge 2$};\\
0,&\text{otherwise}.
\endcases
$$
\endproclaim

Theorem 7.3 suggests the following 

\proclaim{7.4. Problem} Classify $ms(X,r)$ for compact subsets $X$ of
positive Lebesgue measure in $\IR^n$, where $n>2$. In particular, calculate
$ms(X,r)$ for all regular polyhedra in $\IR^n$ for $n>2$. Calculate
$ms(T,2)$ (or equivalently, $ms(S_4\ti I,2)$) for the regular tetrahedron
$T$ in $\IR^3$.
\endproclaim

Observe that the last statement of Theorem 7.3, Theorem 6.2, and known
results on classification of regular polyhedra \cite{20, 12.6} 
imply the following
asymptotic estimation first noticed by Ya.Vorobets.

\proclaim{7.5. Proposition} If $m\ge 2$ and $r=2^k$ for some $k\in\IN$, then
$$
\frac1{r^2}\ge ms(X,r)\ge\frac1{r^2}\frac{2^m-r}{2^m-1}.
$$
for every regular convex polyhedron of dimension $\dim X\ge 2m-1$.

\endproclaim

\heading
8. Invariant $ms([0,1],r)$.
\endheading

In this section we consider one of the most difficult and challanging open
problems related to our subject --- calculating the value of the invariant
$ms([0,1],r)$. To feel why this is important we consider the function
$MS([n],r)$ defined as follows.

For natural $n$ and $r$ let $MS([n],r)$ be the maximal $k\in\IN$ such that
for every $r$-coloring of the discrete interval $[n]=\{1,\dots,n\}$ there
exists a symmetric monochromatic subset $A\subset[n]$ of size $|A|\ge k$.
Recall that under a symmetric subset of the real line we understand a 
subset $A\subset \IR$ such that $A=2x-A$ for some $x\in\IR$. 

Note that the function $MS([n],r)$ resembles the well known Van der
Waerden function $M(n,r)$ assigning to each $n$ and $r$ the maximal $k$
such that for every $r$-coloring of $[n]$ there exists a monochromatic
arithmetic progression of length $k$. Since each arithmetic progression is
a symmetric set, we get $M(n,r)\le MS([n],r)$ for any $n,r\in\IN$. It is
know that for a given $r\ge 2$ the Van der Waerden function $M(n,r)$ very
slowly approaches the infinity as $n\to \infty$. In the contract, the
function $MS([n],r)$ grows like a linear function of $n$, that is, for
every $r\in\IN$ there exists a positive limit
$\lim_{n\to\infty}ms([n],r)$, where $ms([n],r)=\frac{MS([n],r)}{n}$. As
expected, this limit is equal to the invariant $ms([0,1],r)$ defined as
the least upper bound of $\e>0$ such that for every measurable
$r$-coloring of $[0,1]$ there exists a monochromatic symmetric subset
$A\subset[0,1]$ of Lebesgue measure $\lambda(A)\ge\e$.

The exact values of the invariants $ms([0,1],r)$ are not known. However,
for $r=2$ by a true tour de force it was proved the estimations 
$\frac1{4+\sqrt 6}\le ms([0,1],2)<\frac5{24}$, see \cite{7, 6.10}. 
The lower bound was obtained
by a subtle geometric averaging arguments, while for the upper bound
it was necessary to create the machinery of so-called blurred colorings.
To define them we first look at ordinary colorings from a bit different
point of view.

Observe that every $r$-coloring $\chi:X\to[r]$ of a space $X$ can be
identified with the collection $\{\chi_i\}_{i=1}^r$ of functions
$\chi_i:X\to\{0,1\}$ such that $\sum_{i=1}^r\chi_i\equiv1$ (and $\chi_i(x)=1$
iff $\chi(x)=i$ for $x\in X$). Under a {\it (measurable) blurred
$r$-coloring} of $X$ we understand a collection $\{\chi_i\}_{i=1}^r$ of
(measurable) functions $\chi_i:X\to[0,1]$ such that
$\sum_{i=1}^r\chi_i\equiv 1$.

Next, suppose $X$ is a space endowed with a probability measure $\mu$ and
$\Cal S$ is a set of measure-preserving {\it involutive} bijections of
$X$. Observe that in this case for every $r\in\IN$ the invariant
$ms(X,\Cal S,r)$ is equal to
$$
\inf_{\chi}\sup_{s\in\Cal S}\max_{1\le i\le
r}\int_X\chi_i(x)\chi_i(s(x))d\mu, $$
where the infimum is taken over all measurable $r$-colorings
$\chi=\{\chi_i\}_{i=1}^r$ of $X$. Taking this infimum over all measurable
blurred $r$-colorings of $X$, we obtain the definition of the blurred
invariant $bms(X,\Cal S,r)$. Clearly, $bms(X,\Cal S,r)\le ms(X,\Cal S,r)$.
In the sequel we shall omit the symbol $\Cal S$ if it is clear from the
context and shall use the notation $ms(X,r)$ and $bms(X,r)$ in place of
$ms(X,\Cal S,r)$ and $bms(X,\Cal S,r)$, respectively.

In such a way we define the blurred invariants $bms([n],r)$ and
$bms([0,1],r)$ of discrete and continuous intervals (on the discrete
interval $[n]$ we consider the measure $\mu(A)=\frac{|A|}n$, where
$A\subset[n]$). The following theorem proved in \cite{7, \S6} describes
interplay between the blurred and ordinary invariants $ms$.

\proclaim{8.1. Theorem} 
\roster
\item $ms([0,1],r)=\lim_{n\to\infty} ms([n],r)=\inf_{n\in\IN}ms([n],r)$
for any $r\in\IN$;
\item $bms([0,1],r)=\lim_{n\to\infty}bms([n],r)=\inf_{n\in\IN}bms([n],r)$
for any $r\in\IN$;
\item $bms([n],r)\le ms([n],r)$ for any $n,r\in\IN$;
\item $bms([0,1],r)=ms([0,1],r)$ for $r>1$;
\item $\frac1{r^2+r\sqrt{r^2-r}}\le ms([0,1],r)<\frac1{r^2}$ for $r>1$;
\item $\frac1{4+\sqrt 6}\le ms([0,1],2)<\frac5{24}$;
\item $\lim_{r\to\infty}ms([0,1],r)\,r^2=\inf_{r\in\IN}
ms([0,1],r)\,r^2=c$ for some constant $\frac12\le c<\frac5{6}$.
\endroster
\endproclaim

Values $ms([n],2)$ were calculated by computer for all $n<25$. In
contrast, values $bms([n],2)$ are not known even for small $n>3$.
For $n=4$ we know that $bms([4],2)\le \frac5{24}$ --- this explains why
the constant $\frac5{24}$ appeared in Theorem 8.1. Answers (even
particular) to the following problems would be very interesting.

\proclaim{8.2. Problem} Calculate the values $bms([n],r)$ (at least for
small $n$ and $r=2$).
\endproclaim

\proclaim{8.3. Problem} Calculate the value $ms([0,1],r)$ (at least for
$r=2$).
\endproclaim

\proclaim{8.4. Problem} Calculate the constant $c$ from Theorem 8.1(7).
\endproclaim

\proclaim{8.5. Question} How quickly the sequences
$\{ms([n],r)\}_{n=1}^\infty$ and $\{bms([n],r)\}_{n=1}^\infty$ converge to
their limit $ms([0,1],r)$?
\endproclaim

\enddocument

\bye